# ON THE POSTERIOR DISTRIBUTION OF THE NUMBER OF COMPONENTS IN A FINITE MIXTURE

BY AGOSTINO NOBILE

*University of Glasgow*

The posterior distribution of the number of components $k$ in a finite mixture satisfies a set of inequality constraints. The result holds irrespective of the parametric form of the mixture components and under assumptions on the prior distribution weaker than those routinely made in the literature on Bayesian analysis of finite mixtures. The inequality constraints can be used to perform an "internal" consistency check of MCMC estimates of the posterior distribution of $k$ and to provide improved estimates which are required to satisfy the constraints. Bounds on the posterior probability of $k$ components are derived using the constraints. Implications on prior distribution specification and on the adequacy of the posterior distribution of $k$ as a tool for selecting an adequate number of components in the mixture are also explored.

**1. Introduction.** Finite mixture distributions have received much attention in the last decade, as a tool for modeling population heterogeneity and especially as a conceptually simple way of relaxing distributional assumptions. Undoubtedly the development of Markov chain Monte Carlo methods has played an essential catalytic role. A survey of the theory and applications of finite mixtures pre-MCMC is provided by Titterington, Smith and Makov (1985), and a more recent introduction to the topic is Robert (1996). Progress has been particularly evident in the Bayesian approach, where it began with the Gibbs sampling algorithm of Diebolt and Robert (1994) for estimating the parameters of a mixture with a fixed number of components. Subsequent work has considered the number of components $k$ as an object of inference, either using tests to select an adequate number of components or summarizing the uncertainty about it by reporting its posterior distribution. Carlin and Chib (1995) and Raftery (1996) have proposed using Bayes factors to test $k$ against $k + 1$ components and they have described MCMC methods to compute the necessary marginal likelihoods. The paper by Raftery contains a summary of such methods. Mengersen and Robert (1996) also assume a testing perspective, but use the Kullback–Leibler divergence as a measure of distance between models with $k$ and $k + 1$ components. Nobile (1994), Phillips and Smith (1996), Richardson and Green (1997), Roeder and Wasserman (1997) and Stephens (2000) have put a prior distribution on the number of components and obtained MCMC estimates of the posterior. Besides representing uncertainty about $k$, its posterior distribution can also be used to mix models with different numbers of components, leading to model mixing predictions of future observables. Nobile (1994) attempted to estimate the marginal likelihoods of each model separately and then formed an estimate of the posterior of $k$ using Bayes' theorem. Roeder and Wasserman (1997) proposed to approximate the marginal likelihoods using the Schwarz criterion. Although their methods differ considerably, Phillips and Smith (1996), Richardson and Green (1997) and Stephens (2000) share a common approach consisting of running an MCMC sampler on a composite model, with jumps between submodels that allow the sampler to change the number of components in the mixture. Then the posterior of $k$ can be estimated by the relative amount of simulation time spent by the sampler in each submodel.

In this paper I show that, under some conditions on the prior distribution, the marginal likelihoods of finite mixture models with a different number of components satisfy a set of inequality constraints. Besides its theoretical interest, the result provides a means of performing a check of "internal" consistency of MCMC estimates of the marginal likelihoods, or of the marginal likelihoods implicit in MCMC estimates of the posterior of $k$.

**2. The model.** Let $x = \{x_1, \ldots, x_n\}$ be a sequence of (possibly vector-valued) random variables and assume that the $x_i$'s are independent and identically distributed with probability density function (with respect to some underlying measure) given by

$$(1) \qquad\qquad f(x_i) = \sum_{j=1}^{k} \lambda_j p_j(x_i).$$

Model (1) is called a "finite mixture distribution." The mixture weights $\lambda_j$ are the probabilities that the random variable $x_i$ follows any of $k$ alternative distributions, with densities $p_j(\cdot)$, called the "mixture components." In this paper attention is restricted to the case where the number of components $k$, the weights $\lambda_j$ and the components $p_j(\cdot)$ are all unknown. It is assumed, however, that the densities $p_j(\cdot)$ belong to some specified parametric family, allowed to vary with $j$. Thus $p_j(x_i) = p_j(x_i|\theta_j)$, where $\theta_j$ is the vector of parameters of the $j$th mixture component.

As stated, model (1) is somewhat ambiguous, since the meaning of mixture weights and mixture components is completely specified only when $k$ is fixed; for instance, the expression "the weight of the second component" seems to have a different meaning when $k = 2$ than it has when $k = 5$. In order to make explicit the dependence on $k$ of mixture weights and components, rewrite model (1) as follows:

$$f(x_i|k, \lambda_k, \theta_k) = \sum_{j=1}^{k} \lambda_{jk} p_{jk}(x_i|\theta_{jk}), \qquad i = 1, \ldots, n,$$

where $\lambda_k = (\lambda_{1k}, \ldots, \lambda_{kk})^{\top}$ and $\theta_k = (\theta_{1k}, \ldots, \theta_{kk})^{\top}$. On occasion $\lambda = (\lambda_1, \lambda_2, \ldots)^{\top}$ and $\theta = (\theta_1, \theta_2, \ldots)^{\top}$ will be used. In principle this formulation allows the parametric family of the component to change with $j$ and $k$.

Conditional on $k$, let $g_i$ be an integer in $\{1, \ldots, k\}$ denoting the unknown component from which the $i$th observation $x_i$ proceeds. The unobserved vector $g = (g_1, \ldots, g_n)^{\top}$ has been called the "membership vector" or "allocation vector" or "configuration vector" of the mixture. If one conditions on $g$, the distribution of $x_i$ is simply given by the $g_i$th component in the mixture,

$$f(x|k, g, \theta_k) = \prod_{i=1}^{n} p_{g_i,k}(x_i|\theta_{g_i,k}).$$

The complete specification of the Bayesian finite mixture model requires a prior distribution for all the unknown quantities. The prior on $k$, denoted by $\pi(k)$, has support on (a subset of) the positive integers and may involve a hyperparameter. Given $k$, the weights $\lambda_k = (\lambda_{1k}, \ldots, \lambda_{kk})^{\top}$ are assumed to have the $\mathrm{Dir}(\alpha_{1k}, \ldots, \alpha_{kk})$ prior distribution, where the hyperparameters $\alpha_k = (\alpha_{1k}, \ldots, \alpha_{kk})^{\top}$ are positive constants. Although other priors could be used for the weights, the Dirichlet distribution has become a standard choice. The allocations $g_i$ are conditionally independent given $k$ and $\lambda_k$ with $\Pr[g_i = j|k, \lambda_k] = \lambda_{jk}$. Given $k$, independent priors are usually assumed for the component parameters $\theta_{jk}$,

$$\pi(\theta_k|k, \phi_k) = \prod_{j=1}^{k} \pi_{jk}(\theta_{jk}|\phi_{jk}),$$

where $\phi_{jk}$ is the set of hyperparameters in the prior distribution of $\theta_{jk}$ and $\phi_k = (\phi_{1k}, \ldots, \phi_{kk})^{\top}$. In general the components' hyperparameters $\phi_{jk}$ can vary with $k$, so that substantive prior information distinguishing the components and depending on their number $k$ can be accommodated. Similarly, the functional form of the prior $\pi_{jk}(\cdot)$ may change with $j$ and $k$, since the component parametric family may too. Dependence on $k$ is, however, ruled out by the assumptions introduced in Section 3.

In summary, the joint distribution of the data and all unknowns in the model is

$$(2) \qquad \begin{aligned} &f(x, \theta, g, \lambda, k) \\ &= \pi(k)\pi(\lambda_k|k, \alpha_k)f(g|k, \lambda_k)\pi(\theta_k|k, \phi_k)f(x|k, g, \theta_k). \end{aligned}$$

In the sequel, attention is focused on a model obtained by integrating the parameters $\lambda_k$ and $\theta_k$ out of model (2). Integrating the weights out of the model yields

$$f(g|k, \alpha_k) = \int f(g|k, \lambda_k)\pi(\lambda_k|k, \alpha_k)\,d\lambda_k$$

$$= \int \prod_{j=1}^{k} \lambda_{jk}^{n_j} \frac{\Gamma(\alpha_{0k})}{\prod_{j=1}^{k}\Gamma(\alpha_{jk})} \prod_{j=1}^{k} \lambda_{jk}^{\alpha_{jk}-1}\,d\lambda_k$$

$$(3) \qquad = \frac{\Gamma(\alpha_{0k})}{\Gamma(\alpha_{0k}+n)} \prod_{j=1}^{k} \frac{\Gamma(\alpha_{jk}+n_j)}{\Gamma(\alpha_{jk})},$$

where $\alpha_{0k} = \sum_{j=1}^{k}\alpha_{jk}$, $n_j = n_j(g) = \operatorname{card}\{A_j\}$ and $A_j = \{i : g_i = j\}$ is the index set of the observations allocated to the $j$th component. One can also, at least in principle, integrate the component parameters out of the model,

$$f(x|k, g, \phi_k) = \int f(x|k, g, \theta_k)\pi(\theta_k|k, \phi_k)\,d\theta_k$$

$$= \int \prod_{i=1}^{n} p_{g_i,k}(x_i|\theta_{g_i,k}) \prod_{j=1}^{k}\pi_{jk}(\theta_{jk}|\phi_{jk})\,d\theta_k$$

$$(4) \qquad = \prod_{j=1}^{k}\int \prod_{i \in A_j} p_{jk}(x_i|\theta_{jk})\pi_{jk}(\theta_{jk}|\phi_{jk})\,d\theta_{jk}$$

$$(5) \qquad = \prod_{j=1}^{k} q_{jk}(x^j|\phi_{jk}),$$

where $x^j = \{x_i : i \in A_j\}$ comprises the observations that, according to the membership vector $g$, are from the $j$th component and $q_{jk}(x^j|\phi_{jk})$ is a short way of writing the integral in (4), that is, the marginal density of these observations after the parameter $\theta_{jk}$ has been integrated out.

In the end the joint distribution of the data and unknowns is given by

$$(6) \qquad f(x, g, k|\phi, \alpha) = f(x|k, g, \phi_k)f(g|k, \alpha_k)\pi(k),$$

where $\phi = (\phi_1, \phi_2, \dots)^\top$ and $\alpha = (\alpha_1, \alpha_2, \dots)^\top$. Even though the $\phi$'s and $\alpha$'s are fixed constants, I prefer, with a slight abuse of notation, to list them explicitly to the right of the conditioning bars, as it is important to recall that they enter in the expressions in (6). The posterior distribution of the number of components is

$$\pi(k|x, \phi, \alpha) \propto \pi(k) f(x|k, \phi_k, \alpha_k).$$

The marginal likelihoods $f(x|k, \phi_k, \alpha_k)$, which will also be denoted as $f_k$ for short, are given by

$$(7) \quad f_k = f(x|k, \phi_k, \alpha_k) = \sum_{g \in \mathcal{G}_k} f(x|k, g, \phi_k) f(g|k, \alpha_k), \qquad k = 1, 2, \ldots,$$

where the sum extends over the lattice $\mathcal{G}_k = \{g : g_i \in \{1, \ldots, k\}, i = 1, \ldots, n\}$, the set of membership vectors with components at most $k$. Representation (7) demonstrates the great advantage of working with model (6) rather than model (2). Using (7) it becomes possible to compare the contributions of the *same* membership vector $g$ to different $f_k$'s. This leads to linking together the marginal likelihoods and deriving a set of linear inequalities satisfied by them.

**3. Linking the marginal likelihoods.** In this section it is shown that, under certain conditions on the prior distribution, the marginal likelihoods $f_k$ in (7) satisfy a set of constraints. Intuitively, the approach will consist of breaking up the sum over $\mathcal{G}_k$ in (7) into many terms and then showing that some of them can be rewritten as sums over $\mathcal{G}_t$ with $t < k$. The following assumptions will be made throughout.

ASSUMPTION A.1. The Dirichlet hyperparameter of any mixture weight does not change with the number of components:

$$\alpha_{jk} = \alpha_{jj}, \qquad j = 1, \ldots, k - 1, k = 2, 3, \ldots.$$

ASSUMPTION A.2. The properties of any mixture component (parametric family and parameter prior distribution) do not change with the number of components:

$$p_{jk}(\cdot|\cdot) = p_{jj}(\cdot|\cdot), \qquad \pi_{jk}(\cdot|\cdot) = \pi_{jj}(\cdot|\cdot), \qquad \phi_{jk} = \phi_{jj},$$

$$j = 1, \ldots, k - 1, k = 2, 3, \ldots.$$

The assumptions impose a coherency requirement. Not only the $j$th component "remains the same" whether there are $k$ or $k' < k$ components in the mixture (Assumption A.2), but the probability distribution of the ratio between the weight of the $j$th component and the sum of the weights of the first $k'$ components also remains unchanged (Assumption A.1). Because of Assumptions A.1 and A.2, when referring to a certain component one can do so without specifying the number of components in the mixture.

Begin by noticing that the space of membership vectors $\mathcal{G}_k$ in (7) can be partitioned as follows:

$$(8) \qquad \mathcal{G}_k = \bigcup_{t=1}^{k} \mathcal{G}_t^{\star}, \qquad \mathcal{G}_t^{\star} \cap \mathcal{G}_s^{\star} = \varnothing, \qquad t \neq s,$$

where $\mathcal{G}_t^{\star}$ is the set of membership vectors that assigns at least one observation to the $t$th component and none to higher components: $\mathcal{G}_t^{\star} = \{g \in \mathcal{G}_t : \exists\, i \text{ s.t. } g_i = t\}$.

DEFINITION 3.1. Let $f_t^{\star}$ be the portion of $f_t$ that accounts for the membership vectors $g$ that allocate at least one observation to component $t$ and none to higher components (components lower than $t$ may be empty),

$$(9) \qquad f_t^{\star} = \sum_{g \in \mathcal{G}_t^{\star}} f(x|t,g,\phi_t) f(g|t,\alpha_t).$$

Clearly $f_1^{\star} = f_1$.

In the sequel use will be made of the following conditions.

CONDITION C.1. For all $g \in \mathcal{G}_t^{\star}$ with $t < k$,

$$f(x|k,g,\phi_k) = f(x|t,g,\phi_t).$$

CONDITION C.2. For all $g \in \mathcal{G}_t^{\star}$ with $t < k$,

$$\frac{f(g|k,\alpha_k)}{f(g|t,\alpha_t)} = a_{kt} \qquad \text{constant.}$$

LEMMA 3.1. *Under Assumptions* A.1 *and* A.2, *the model of Section* 2 *satisfies Conditions* C.1 *and* C.2 *with*

$$(10) \qquad a_{kt} = \frac{\Gamma(\alpha_{0k})}{\Gamma(\alpha_{0k}+n)} \frac{\Gamma(\alpha_{0t}+n)}{\Gamma(\alpha_{0t})}.$$

PROOF. To verify Condition C.1, recall (5): $f(x|k,g,\phi_k) = \prod_{j=1}^{k} q_{jk}(x^j|\phi_{jk})$. All $g \in \mathcal{G}_t^{\star}$, $t < k$, allocate no observations to components larger than the $t$th one: $x^j = \varnothing, j > t$. Therefore the product in (5) extends from 1 to $t$ only. Moreover, Assumption A.2 implies that, for $j \in \{1,\ldots,t\}$, $q_{jk}(\cdot|\cdot) = q_{jt}(\cdot|\cdot)$ and $\phi_{jk} = \phi_{jt}$. Hence $f(x|k,g,\phi_k) = \prod_{j=1}^{t} q_{jt}(x^j|\phi_{jt}) = f(x|t,g,\phi_t)$. As for Condition C.2, from (3) one has

$$\frac{f(g|k,\alpha_k)}{f(g|t,\alpha_t)} = \frac{\Gamma(\alpha_{0k})}{\Gamma(\alpha_{0k}+n)} \prod_{j=1}^{k} \frac{\Gamma(\alpha_{jk}+n_j)}{\Gamma(\alpha_{jk})} \bigg/ \frac{\Gamma(\alpha_{0t})}{\Gamma(\alpha_{0t}+n)} \prod_{j=1}^{t} \frac{\Gamma(\alpha_{jt}+n_j)}{\Gamma(\alpha_{jt})}.$$

Again, for all $g \in \mathcal{G}_t^\star$ and $j > t$, $A_j = \varnothing$ so that $n_j = 0$. Hence the last $k - t$ terms in the product in the numerator are 1. Also, from Assumption A.1, $\alpha_{jk} = \alpha_{jt}$, $j = 1, \ldots, t$. Therefore C.2 holds with $a_{kt}$ given by (10). $\square$

The following result may be considered as an appetizer.

THEOREM 3.1. *Let $f_k$ and $f_t^\star$ be as in* (7) *and* (9) *and assume that Conditions* C.1 *and* C.2 *hold. Then*

$$(11) \qquad f_k = \sum_{t=1}^{k} a_{kt} f_t^\star.$$

*Moreover,*

$$(12) \qquad f_k = a_{k,k-1} f_{k-1} + f_k^\star.$$

PROOF. Equation (7) can be rewritten as $f_k = \sum_{t=1}^{k} \sum_{g \in \mathcal{G}_t^\star} f(x|k, g, \phi_k) \times f(g|k, \alpha_k)$ because of the partition of $\mathcal{G}_k$ in (8). Now use Conditions C.1 and C.2 and the definition of $f_t^\star$ in (9) to obtain (11). A little more algebra yields (12):

$$f_k = \sum_{t=1}^{k} a_{kt} f_t^\star = f_k^\star + \sum_{t=1}^{k-1} a_{kt} \frac{f(g|k-1, \alpha_{k-1})}{f(g|k-1, \alpha_{k-1})} f_t^\star$$

$$= f_k^\star + \sum_{t=1}^{k-1} a_{k-1,t} \frac{f(g|k, \alpha_k)}{f(g|k-1, \alpha_{k-1})} f_t^\star$$

$$= f_k^\star + a_{k,k-1} \sum_{t=1}^{k-1} a_{k-1,t} f_t^\star = a_{k,k-1} f_{k-1} + f_k^\star. \qquad \square$$

Theorem 3.1 provides two representations of $f_k$. In (11) it is given as a linear combination of the "no empty last component" portions of the marginal likelihoods of models with $1, 2, \ldots, k$ components. In (12) it is written as the "no empty last component" portion of the marginal likelihood of the $k$-components model plus a fraction of the marginal likelihood of the model with one fewer component. Much of the remainder of this section is devoted to deriving a result stronger than Theorem 3.1. This is achieved by exploiting additional symmetry left as yet untapped; some mixture components may have identical characteristics. The first step consists in grouping the mixture components into classes of "alike" components.

DEFINITION 3.2. Say that two mixture components $j$ and $k$ are alike or equivalent if $\alpha_{jj} = \alpha_{kk}$, $p_{jj}(\cdot|\cdot) = p_{kk}(\cdot|\cdot)$, $\pi_{jj}(\cdot|\cdot) = \pi_{kk}(\cdot|\cdot)$ and $\phi_{jj} = \phi_{kk}$.

The above definition induces a partition of the components into classes of equivalence, with two components being in the same class if they are alike. It may help intuition to regard the observations as balls being placed in a sequence of colored boxes, with boxes of the same color being equivalent. Let $\mathcal{C}(m)$ be the $m$th equivalence class and let $m.h$ be the index of the $h$th smallest component in $\mathcal{C}(m)$. The classes are ordered so that $\mathcal{C}(m)$ precedes $\mathcal{C}(r)$ if $m.1 < r.1$. Each class contains either a finite number of components, possibly one, or countably many components, possibly all. Let $N(t)$ be the number of equivalence classes formed by components 1 through $t$. Also, let $i(t)$ be the index of the equivalence class to which component $t$ belongs, so that $\mathcal{C}(i(t))$ is the class of components that are equivalent to component $t$. Finally, let $c(m,t)$ be the number of components in $\mathcal{C}(m)$ that are no larger than $t$ and let $c(m)$ be its total number of components: $c(m) = \sup_t c(m,t)$. One extreme case often considered in the literature is that of just one equivalence class: there is no prior information distinguishing the components. In this case $N(t) \equiv 1$, $\mathcal{C}(1) = \{1, 2, \dots\}$, $1.h = h$, $i(t) \equiv 1$, $c(1,t) = t$ and $c(1) = \infty$. The other extreme case arises when each class contains only one component; $N(t) = t$, $\mathcal{C}(m) = \{m\}$, $m.1 = m$, $i(t) = t$, $c(m,t) = I(m \le t)$ and $c(m) \equiv 1$, with $I(\cdot)$ the indicator function.

DEFINITION 3.3. For any membership vector $g \in \mathcal{G}_t^\star$, define its class occupancy pattern as the vector $h = (h_1, \dots, h_{N(t)})^\top$, where $h_m$ is the number of nonempty components in class $\mathcal{C}(m)$.

Let $H_t: \mathcal{G}_t^\star \longrightarrow \{0, 1, 2, \dots\}^{N(t)}$ be the mapping which associates to each $g \in \mathcal{G}_t^\star$ its class occupancy pattern $h$. Since the domain of $H_t$ is $\mathcal{G}_t^\star$, component $t$ is nonempty, hence $h_{i(t)} \ge 1$; also, the number of nonempty components cannot exceed the number of observations. Therefore, the range of the mapping, $\mathcal{H}_t = H_t(\mathcal{G}_t^\star)$, consists of the $N(t)$-dimensional vectors $h$ satisfying

$$(13) \qquad \sum_{m=1}^{N(t)} h_m \le n, \qquad h_m \in \begin{cases} \{1, 2, \dots, c(m,t)\}, & \text{if } m = i(t), \\ \{0, 1, \dots, c(m,t)\}, & \text{otherwise.} \end{cases}$$

If $\sum_{m=1}^{N(t)} h_m < t$, some mixture components in $\{1, \dots, t\}$ are empty. This suggests that it may be possible to accommodate the class occupancy pattern $h$ using fewer than $t$ components.

DEFINITION 3.4. For any class occupancy $h$, let $s = s(h)$ be the smallest integer such that the mixture components from 1 to $s$ comprise at least $h_m$ components in $\mathcal{C}(m)$, $m = 1, \dots, \text{card}(h)$,

$$(14) \qquad s = s(h) = \min\{r : c(m,r) \ge h_m, m = 1, \dots, \text{card}(h)\},$$

where $\text{card}(h)$ is the number of elements of $h$. If $h \in \mathcal{H}_t$, then $\text{card}(h) = N(t)$ and $s \le t$.

The symbol $s$ will be exclusively used to denote the function defined in (14). For any $h \in \mathcal{H}_t$, $s$ is the smallest number of components needed to accommodate $h$, so that $h \in \mathcal{H}_s$ too, under the convention that trailing 0's in $h$ are dropped. For instance, suppose that $t = 6$, $\mathcal{C}(1) \supset \{1, 2, 3, 6\}$, $\mathcal{C}(2) \supset \{4\}$, $\mathcal{C}(3) \supset \{5\}$, so that $N(6) = 3$. If $h = (2, 1, 0)^\top$ then only three components are nonempty and $s = 4$. Dropping the trailing 0 in $h$, $h = (2, 1)^\top \in \mathcal{H}_4$.

DEFINITION 3.5. Let $\mathcal{H}_r^t = \{h \in \mathcal{H}_t : r = s(h)\}$ be the (possibly empty) subset of class occupancies $\mathcal{H}_t$ which can be accommodated with $r \leq t$ components.

The set of class occupancies of the membership vectors in $\mathcal{G}_t^\star$ can be partitioned as follows:

$$(15) \qquad \mathcal{H}_t = \bigcup_{r=1}^{t} \mathcal{H}_r^t, \qquad \mathcal{H}_r^t \cap \mathcal{H}_q^t = \varnothing, \qquad r \neq q.$$

If $h \in \mathcal{H}_r^t$, then $s(h) = r$ so that $h \in \mathcal{H}_r$ too, and hence $h \in \mathcal{H}_r^r$. This shows that

$$(16) \qquad \mathcal{H}_r^t \subset \mathcal{H}_r^r, \qquad r < t.$$

DEFINITION 3.6. Let $\mathcal{G}_h^t$ with $t \geq s(h)$ be the subset of $\mathcal{G}_t^\star$ consisting of membership vectors with class occupancy pattern $h : \mathcal{G}_h^t = H_t^{-1}(h)$.

Clearly, $\{\mathcal{G}_h^t, h \in \mathcal{H}_t\}$ is a partition of $\mathcal{G}_t^\star$:

$$(17) \qquad \mathcal{G}_t^\star = \bigcup_{h \in \mathcal{H}_t} \mathcal{G}_h^t, \qquad \mathcal{G}_h^t \cap \mathcal{G}_v^t = \varnothing, \qquad h \neq v.$$

Consider next the mapping $M_t : \mathcal{G}_h^t \longrightarrow \mathcal{G}_s^\star$ which removes any gap in the sequence of nonempty components within each equivalence class. More precisely, given $g \in \mathcal{G}_h^t$, let $j_{m1} < \cdots < j_{m,h_m}$ be the corresponding nonempty components in $\mathcal{C}(m)$, $m = 1, \ldots, N(t)$. The mapping $M_t$ changes, for all $m$, the components $j_{m1}, \ldots, j_{m,h_m}$ into $m.1, \ldots, m.h_m$, respectively. Denote the range of $M_t$ by $\mathcal{E}_h = M_t(\mathcal{G}_h^t)$, noting that from the definition of $\mathcal{G}_h^t$ it is immediate that $M_t(\mathcal{G}_h^t) = M_r(\mathcal{G}_h^r)$ for any $t, r \geq s(h)$. The mapping $M_t$ does not affect the class occupancy of a membership vector; thus $H_s(\mathcal{E}_h) = \{h\}$, although in general $\mathcal{E}_h$ is a subset of $H_s^{-1}(h) = \mathcal{G}_h^s$. Because of the equivalence of components within each class, the mapping $M_t$ leaves unchanged $f(x, g | t, \phi_t, \alpha_t)$,

$$(18) \qquad f(x | t, g, \phi_t) f(g | t, \alpha_t) = f(x | t, \tilde{g}, \phi_t) f(\tilde{g} | t, \alpha_t), \qquad g \in \mathcal{G}_h^t, \tilde{g} = M_t(g).$$

DEFINITION 3.7. *Let $\gamma_h^t$ be defined as follows:*

$$(19) \qquad \gamma_h^t = \begin{cases} \dfrac{h_{i(t)}}{c(i(t),t)} \displaystyle\prod_{m=1}^{N(t)} \binom{c(m,t)}{h_m}, & h \in \mathcal{H}_t, \\ 0, & h \notin \mathcal{H}_t. \end{cases}$$

LEMMA 3.2. *Any element of $\mathcal{E}_h$ is the image under $M_t$ of $\gamma_h^t$ membership vectors in $\mathcal{G}_h^t$.*

Lemma 3.2 says that $\mathcal{G}_h^t$ consists of $\gamma_h^t$ subsets alike to $\mathcal{E}_h$, except for which $h_m$ components in each class $\mathcal{C}(m)$ are nonempty. Coupled with (18), Lemma 3.2 gives

$$(20) \qquad \sum_{g \in \mathcal{G}_h^t} f(x|t,g,\phi_t) f(g|t,\alpha_t) = \gamma_h^t \sum_{g \in \mathcal{E}_h} f(x|t,g,\phi_t) f(g|t,\alpha_t).$$

DEFINITION 3.8. *Let $f_h^\dagger$ be the portion of $f_s$, $s = s(h)$, that accounts for the membership vectors in $\mathcal{G}_h^s$,*

$$(21) \qquad f_h^\dagger = \sum_{g \in \mathcal{G}_h^s} f(x|s,g,\phi_s) f(g|s,\alpha_s), \qquad s = s(h).$$

The following lemma is instrumental in proving the main result, Theorem 3.2.

LEMMA 3.3. *The function $f_t^\star$ defined in (9) can be rewritten as follows:*

$$(22) \qquad f_t^\star = \sum_{r=1}^{t} a_{tr} \sum_{h \in \mathcal{H}_r^r} \frac{\gamma_h^t}{\gamma_h^r} f_h^\dagger.$$

THEOREM 3.2. *Suppose that Conditions C.1 and C.2 are verified. Then*

$$(23) \qquad f_k = \sum_{r=1}^{k} a_{kr} \sum_{h \in \mathcal{H}_r^r} \gamma_h^{r,k} f_h^\dagger,$$

*where $f_h^\dagger$ is defined in (21),*

$$(24) \qquad \gamma_h^{r,k} = \frac{1}{\gamma_h^r} \sum_{t=r}^{k} \gamma_h^t$$

*and $\gamma_h^t$ is given in (19). Moreover,*

$$(25) \qquad f_k = a_{k,k-1} f_{k-1} + \sum_{r=1}^{k} a_{kr} \sum_{h \in \mathcal{H}_r^r} \frac{\gamma_h^k}{\gamma_h^r} f_h^\dagger.$$

It is worthwhile to consider explicitly the cases where all components are equivalent and where no two components are equivalent.

PROPOSITION 3.1. *Suppose that Conditions* C.1 *and* C.2 *are satisfied and that all mixture components are equivalent. Then*

$$(26) \qquad f_k = \sum_{h=1}^{k \wedge n} \binom{k}{h} a_{kh} f_h^\dagger$$

$$(27) \qquad = a_{k,k-1} f_{k-1} + \sum_{h=1}^{k \wedge n} \binom{k-1}{h-1} a_{kh} f_h^\dagger.$$

PROOF. Recall that if all the components are equivalent then $N(t) \equiv 1$, $c(1, t) = t$ and $i(t) \equiv 1$. Therefore the class occupancy $h$ is a scalar, the number of nonempty components in the unique equivalence class. From formula (13) the range of $h$ is $\mathcal{H}_t = \{1, \ldots, t \wedge n\}$, with $t \wedge n = \min(t, n)$. From Definition 3.4 the smallest number of components needed to accommodate $h$ is $s(h) = h$. Hence $\mathcal{H}_r^t = \{r\}$, $r \leq t \wedge n$, and $\mathcal{H}_r^t = \varnothing$, $r > t \wedge n$. Here the range of $M_t$ is $\mathcal{E}_h = \mathcal{G}_h^h$, the subset of $\mathcal{G}_h$ consisting of membership vectors that allocate at least one observation to each component, while (21) gives the part of the marginal likelihood $f_h$ corresponding to no empty components,

$$(28) \qquad f_h^\dagger = \sum_{g \in \mathcal{G}_h^h} f(x|h, g, \phi_h) f(g|h, \alpha_h).$$

In this case expression (23) becomes $f_k = \sum_{h=1}^{k \wedge n} a_{kh} \gamma_h^{h,k} f_h^\dagger$. From (19) one has $\gamma_h^t = \binom{t-1}{h-1}$ so that $\gamma_h^{h,k} = \sum_{t=h}^{k} \binom{t-1}{h-1} = \binom{k}{h}$ and (26) follows. Equation (27) can be derived from (25) after making substitutions similar to the ones performed to obtain (26). □

Formula (26) provides a representation of the marginal likelihood of $k$ components as a linear combination of the portions of marginal likelihoods corresponding to no empty components.

PROPOSITION 3.2. *Suppose that Conditions* C.1 *and* C.2 *hold and that no two mixture components are equivalent. Then*

$$(29) \qquad f_k = \sum_{t=1}^{k} a_{kt} f_t^\star$$

$$(30) \qquad = a_{k,k-1} f_{k-1} + f_k^\star.$$

The proof is left as an exercise for the interested reader.

Note that the conclusion of Proposition 3.2 coincides with that of Theorem 3.1, if no two components are equivalent there is no additional symmetry to be exploited beyond what is assumed by Theorem 3.1. The following corollary summarizes some special cases.

COROLLARY 3.1. *For the model of Section 2, under Assumptions* A.1 *and* A.2, *one has the following*:

(i) *representations* (23) *and* (25) *hold with* $a_{kr}$ *as given in* (10);

(ii) *in the special case where all mixture components are equivalent with the Dirichlet prior on the mixture weights having hyperparameter* $\alpha_{jk} = \alpha$, *one has*

$$(31) \quad f_k = \sum_{h=1}^{k \wedge n} \binom{k}{h} \frac{\Gamma(k\alpha)}{\Gamma(k\alpha + n)} \frac{\Gamma(h\alpha + n)}{\Gamma(h\alpha)} f_h^\dagger$$

$$(32) \quad = \prod_{i=1}^{n} \left( \frac{k\alpha - \alpha - 1 + i}{k\alpha - 1 + i} \right) f_{k-1} + \sum_{h=1}^{k \wedge n} \binom{k-1}{h-1} \frac{\Gamma(k\alpha)}{\Gamma(k\alpha + n)} \frac{\Gamma(h\alpha + n)}{\Gamma(h\alpha)} f_h^\dagger;$$

(iii) *in case* (ii) *above with* $\alpha = 1$ *one has*

$$f_k = \sum_{h=1}^{k \wedge n} \frac{k!}{h!(k-h)!} \frac{(k-1)!}{(k-1+n)!} \frac{(h-1+n)!}{(h-1)!} f_h^\dagger$$

$$= \frac{k-1}{k+n-1} f_{k-1} + \sum_{h=1}^{k \wedge n} \frac{(k-1)!}{(h-1)!(k-h)!} \frac{(k-1)!}{(k-1+n)!} \frac{(h-1+n)!}{(h-1)!} f_h^\dagger.$$

The representations of the marginal likelihoods $f_k$ provided in Theorems 3.1 and 3.2 and its corollaries lead to a set of linear constraints on the $f_k$'s. Solving the triangular system (11) for the $f_t^{\star}$'s in terms of the $f_k$'s, one obtains (12) $f_k^{\star} = f_k - a_{k,k-1} f_{k-1}$. As the $f_k^{\star}$'s are, from equation (9), sums of strictly positive terms, this implies that

$$(33) \qquad\qquad\qquad f_k > a_{k,k-1} f_{k-1}.$$

The constraints (33) hold no matter how the mixture components partition into classes of equivalence. In the case of no equivalent components treated in Proposition 3.2, the constraints (33) cannot be made any stronger, since by how much $f_k$ exceeds $a_{k,k-1} f_{k-1}$, that is, $f_k^{\star}$, depends on vectors which allocate at least one observation to component $k$. At the opposite extreme of all equivalent components, dealt with in Proposition 3.1, stronger constraints are obtained by solving the triangular system (26) for the $f_h^{\dagger}$'s in terms of the $f_k$'s, and then setting the solution to be positive. These constraints, explicitly derived in formula (36), are stronger than (33) because, of all the $f_h^{\dagger}$'s in the sum $\sum_{h=1}^{k \wedge n}$ in (27), only $f_k^\dagger$ involves vectors allocating observations to the $k$th component. As a very special case, consider equation (26) with $k > n$. Then $f_k$ is a linear combination of $f_1^\dagger, \ldots, f_n^\dagger$. However, $f_1^\dagger, \ldots, f_n^\dagger$ can be obtained by solving (26) with $k = 1, \ldots, n$. Therefore, $f_k$ with $k > n$ is completely determined by the marginal likelihoods $f_1, \ldots, f_n$; this is a much stronger result than is obtainable when no components are equivalent. The general case where only some components are equivalent is covered by Theorem 3.2. As usual the constraints (33) hold, but, contrary to the case of all equivalent components, one cannot solve system (23) for the $f_h^\dagger$'s. Nevertheless, there might be a function of the $f_h^\dagger$'s, finer than $f_t^\star$, such that system (23) can be solved for it.

The remainder of this section deals exclusively with the case where all mixture components are equivalent. The triangular system (26) with $k = 1, \ldots, n$ can be rewritten as

$$(34) \qquad f_k = f_k^\dagger + \sum_{t=1}^{k-1} b_{kt} f_t^\dagger, \qquad k \leq n,$$

with $b_{kt} = \binom{k}{t} a_{kt}$. Denote by $B_n$ the matrix of coefficients of system (34). In this case one can provide a simple explicit expression for the elements of $B_n^{-1}$. The following lemma is needed.

Lemma 3.4. *Consider the $q$-dimensional unit lower triangular matrix $B = \{b_{kt}\}$ with $b_{kt} = \binom{k}{t} a_{kt}$ and $a_{kt}$ as in Condition C.2. Let $C$ be the unit lower triangular matrix with generic element $c_{kt} = (-1)^{k+t} b_{kt}$. Then $B^{-1} = C$.*

Proposition 3.3. *Suppose that Conditions C.1 and C.2 are satisfied and that all mixture components are equivalent. Then*

$$(35) \qquad f_k^\dagger = f_k + \sum_{t=1}^{k-1} (-1)^{k+t} \binom{k}{t} a_{kt} f_t, \qquad k \leq n.$$

Proof. The matrix $B_n$ of the coefficients of system (34) is unit lower triangular with generic element $b_{kt} = \binom{k}{t} a_{kt}$. From Lemma 3.4, the inverse $B_n^{-1}$ has generic element $b^{kt} = (-1)^{k+t} \binom{k}{t} a_{kt}$, $k > t$, and the result follows. $\square$

The following corollary follows immediately from Proposition 3.3 and summarizes some special cases.

Corollary 3.2. *For the model of Section 2, under Assumptions A.1 and A.2, one has the following:*

(i) *if all components are equivalent with Dirichlet prior on the weights having hyperparameter $\alpha$, then*

$$f_k^\dagger = \sum_{t=1}^k (-1)^{k+t} \binom{k}{t} \frac{\Gamma(k\alpha)}{\Gamma(k\alpha+n)} \frac{\Gamma(t\alpha+n)}{\Gamma(t\alpha)} f_t, \qquad k \le n;$$

(ii) *in case* (i) *above with $\alpha = 1$,*

$$f_k^\dagger = \sum_{t=1}^k (-1)^{k+t} \frac{k!}{t!(k-t)!} \frac{(k-1)!}{(k-1+n)!} \frac{(t-1+n)!}{(t-1)!} f_t, \qquad k \le n.$$

Briefly returning to the topic of the inequality constraints on the $f_k$'s, from (35) one has

$$(36) \qquad f_k > \sum_{t=1}^{k-1} (-1)^{k+t+1} \binom{k}{t} a_{kt} f_t, \qquad k \le n.$$

The following section discusses possible uses of these constraints; the present one concludes by addressing the problem of expressing $f_k$ with $k > n$ in Proposition 3.1 in terms of $f_1, \ldots, f_n$.

PROPOSITION 3.4. *Suppose that Conditions* C.1 *and* C.2 *are satisfied and that all mixture components are equivalent. Then*

$$(37) \qquad f_k = \sum_{t=1}^n (-1)^{n-t} \binom{k}{t} \binom{k-t-1}{n-t} a_{kt} f_t, \qquad k > n.$$

**4. Applications.** This section explores some uses of the representations of the marginal likelihoods derived in Section 3.

1. When all mixture components are equivalent, a proper prior on the number of components is necessary in order to have a proper posterior.
2. Bounds on the posterior probability of $k$ mixture components can be derived that hold for any sample of given size and for any family of component distributions.
3. An "internal" consistency check of Markov chain Monte Carlo estimates of the marginal likelihoods $f(x|k)$ can be performed by verifying that they satisfy the constraints. Estimates that fail the check can seemingly be improved by modifying them so that the constraints are satisfied.
4. Expressions can be obtained for the prior and posterior distributions of the number of nonempty components in the mixture, that is, the number of components to which observations are allocated.

Throughout this section attention is focused on the case where all mixture components are equivalent, for a variety of reasons: it is important in practice, it is amenable to a notationally simpler treatment and it leads to stronger results. In order to lighten the notation, the explicit indication of the hyperparameters is abandoned in most of this section, so, for instance, I will write $\pi(k|x)$ and $f(x|k)$ in place of $\pi(k|x, \phi, \alpha)$ and $f(x|k, \phi_k, \alpha_k)$. Fortran and S-PLUS programs used for the computations in this section are available from the author upon request.

4.1. *Proper posterior of $k$.* From Bayes' theorem, the posterior distribution of the number of components is

$$(38) \qquad \pi(k|x) = \frac{\pi(k)f(x|k)}{\sum_{j=1}^{\infty} \pi(j)f(x|j)} = \frac{\pi(k)\sum_{h=1}^{k \wedge n} \binom{k}{h} a_{kh} f_h^{\dagger}}{\sum_{j=1}^{\infty} \pi(j)\sum_{h=1}^{j \wedge n} \binom{j}{h} a_{jh} f_h^{\dagger}},$$

where the representation of the marginal likelihoods given in (26) was used. Since the series in the denominator of (38) is of positive terms, one can change the order of summation to obtain

$$(39) \qquad \pi(k|x) = \frac{\sum_{h=1}^{k \wedge n} f_h^{\dagger} \pi(k) \binom{k}{h} a_{kh}}{\sum_{h=1}^{n} f_h^{\dagger} \{\sum_{j=h}^{\infty} \pi(j) \binom{j}{h} a_{jh}\}}.$$

A proper prior distribution $\pi(k)$ ensures that the posterior is also a proper probability distribution. The following theorem shows that, when all mixture components are equivalent, this condition is not only sufficient but also necessary.

THEOREM 4.1. *Consider the model of Section 2, under Assumptions* A.1 *and* A.2, *and suppose that all mixture components are equivalent. Then the posterior $\pi(k|x)$ of the number of components is a proper distribution if and only if the prior $\pi(k)$ is proper.*

PROOF. The posterior $\pi(k|x)$ is proper if and only if the series in braces in the denominator of (39) converges. Using formula (10) for $a_{jh}$ the series become

$$\frac{\Gamma(\alpha_{0h} + n)}{\Gamma(\alpha_{0h})h!} \sum_{j=h}^{\infty} \frac{j!}{(j-h)!} \frac{\Gamma(\alpha_{0j})}{\Gamma(\alpha_{0j} + n)} \pi(j), \qquad h = 1, \dots, n.$$

Clearly, if the above series converges when $h = n$, it also converges for $h < n$. Thus $\pi(k|x)$ is proper if and only if the following series converges:

$$(40) \qquad \sum_{j=n}^{\infty} \frac{j}{\alpha_{0j} + n - 1} \frac{j-1}{\alpha_{0j} + n - 2} \cdots \frac{j-n+1}{\alpha_{0j}} \pi(j).$$

Since all components are equivalent, $\alpha_{0j} = j\alpha$ for some $\alpha > 0$. Letting $c_j$ denote the generic term of series (40), it is easy to see that

$$(41) \qquad \frac{\pi(j)}{(n\alpha + n - 1)^n} < c_j < \frac{\pi(j)}{\alpha^n}.$$

To prove the right-hand side inequality note that each of the $n$ terms in the product in (40) is smaller than $1/\alpha$. For the left-hand side inequality, note that each of those terms is larger than $(j - n + 1)/(j\alpha + n - 1)$, which in turn is no smaller than $1/(n\alpha + n - 1)$. If the prior on $k$ is proper, it then follows, from the right-hand side inequality of (41) and the comparison test for series, that the posterior is also proper. Similarly, if the prior is not proper, the posterior is also seen to be improper, by an application of the comparison test to the left-hand side inequality of (41). $\square$

4.2. *Bounds on the posterior of $k$.* In this subsection it is assumed that the prior distribution on the number of components is proper. A bound on $\pi(k|x)$ results from the maximization of the right-hand side of (39) with respect to $\{f_h^\dagger\}_{h=1}^n$ subject to $f_h^\dagger \geq 0$. The following result simplifies computations.

PROPOSITION 4.1. *Among the vectors that maximize the right-hand side of (39) there is at least one vector $\{f_h^\dagger\}_{h=1}^n$ with only one nonzero component $f_t^\dagger$, with $t \in \{1, \ldots, k \wedge n\}$.*

Note that the nonzero component of the maximizer in Proposition 4.1 need not be the $(k \wedge n)$th. Also, note that, as a function of $\{f_h^\dagger\}_{h=1}^n$, the right-hand side of (39) is constant over lines through the origin; that is, it is homogeneous of zero degree, so that in computing it one can set $f_t^\dagger = 1$. Proposition 4.1 restricts the set of vectors $\{f_h^\dagger\}_{h=1}^n$ one has to compute to find a maximizer of (39) to the $k \wedge n$ vectors with all but one component equal to 0; one can simply compute the right-hand side of (39) for each of them and then pick the one that yields the maximum value.

The bound thus obtained holds, whatever the distributional form of the components in the mixture, as long as they are all equivalent. Moreover, it only depends on the data through the sample size $n$. As an example, consider the posterior of $k$ for a sample of size $n = 82$, with a discrete uniform prior on $k$ over $\{1, \ldots, k_{\max} = 30\}$ and $\alpha_{jk} = \alpha = 1$ for all $j, k$. A maximizer of (39) with $k = 3$ is $f_3^\dagger = 1$, $f_h^\dagger = 0$, $h \neq 3$. The posterior of $k$ corresponding to the maximizer is reported in Table 1. The bound is $\pi(3|x) \leq 0.8623$. These numerical results remain essentially unchanged for any discrete uniform prior with $k_{\max} \geq 10$.

TABLE 1
*Posterior distribution of $k$ which gives maximum probability to $k = 3$, assuming that $n = 82$, $\pi(k) = k_{\max}^{-1}, k = 1, \ldots, k_{\max} = 30$ and $\alpha = 1$*

| $k$ | 1 | 2 | 3 | 4 | 5 | 6 | 7 | 8 | 9 |
|---|---|---|---|---|---|---|---|---|---|
| $\pi(k|x)$ | 0 | 0 | 0.8623 | 0.1217 | $1.42 \times 10^{-2}$ | $1.63 \times 10^{-3}$ | $1.94 \times 10^{-4}$ | $2.44 \times 10^{-5}$ | $3.26 \times 10^{-6}$ |

Table 2 contains bounds on $\pi(k|x)$ for several values of $k$ and $n$, under a uniform prior on $k$ over $\{1, \ldots, k_{\max} = 50\}$ and $\alpha = 1$. Tables 3 and 4 contain bounds when $\alpha = 2$ and $\alpha = 0.5$, respectively.

Tables 2–4 are still correct, at the reported precision, for any discrete uniform prior on $k$ with $k_{\max} > 50$. Since the bounds involve the data only through the sample size $n$, they provide a glimpse of the strength of the prior distribution. Thus, it is to be expected that, for fixed $k$, the bounds become weaker as sample size increases. Perhaps less obvious is that, for fixed sample size, the bounds become stronger as $k$ increases. An intuitive explanation is as follows. Suppose that the model with $k$ components has considerable posterior mass. The posterior mass of the model with $k + 1$ components is at least in part due to the $k + 1$ copies of $\mathcal{G}_k$ embedded in $\mathcal{G}_{k+1}$, all corresponding to at least one empty component. How large this part is depends on the prior distribution, but it may well increase with $k$ since the larger space contains $k + 1$ copies of the smaller one. The values of the hyperparameters $\alpha_{jk} = \alpha$ also greatly affect the bounds, as one can see by comparing Tables 2–4. Increasing $\alpha$ leads to Dirichlet distributions that make very small mixture weights less probable. In turn this reduces the probability mass assigned by the prior on $g$ to membership vectors with empty components. The effect is to "loosen" the link between the marginal likelihoods of different numbers of components, thus making the bounds weaker. Therefore, a more informative prior on the mixture weights leads to weaker constraints on the posterior of $k$.

TABLE 2
*Bounds on $\pi(k|x)$ for several sample sizes $n$, $\pi(k) = k_{\max}^{-1}, k = 1, \ldots, k_{\max} = 50$, $\alpha = 1$*

| $n$ | $k$ | | | | | | | | | |
|---|---|---|---|---|---|---|---|---|---|---|
| | 1 | 2 | 3 | 4 | 5 | 6 | 7 | 8 | 9 | 10 |
| 20 | 0.9000 | 0.7286 | 0.5299 | 0.3456 | 0.2880 | 0.2419 | 0.1954 | 0.1756 | 0.1505 | 0.1335 |
| 50 | 0.9600 | 0.8847 | 0.7826 | 0.6645 | 0.5414 | 0.4233 | 0.3175 | 0.3119 | 0.2835 | 0.2402 |
| 100 | 0.9800 | 0.9412 | 0.8858 | 0.8170 | 0.7385 | 0.6541 | 0.5677 | 0.4828 | 0.4023 | 0.3322 |
| 500 | 0.9960 | 0.9880 | 0.9762 | 0.9607 | 0.9417 | 0.9193 | 0.8938 | 0.8656 | 0.8350 | 0.8022 |

Table 3
*Bounds on $\pi(k|x)$ for several sample sizes $n$, $\pi(k) = k_{\max}^{-1}, k = 1, \ldots, k_{\max} = 50$, $\alpha = 2$*

| $n$ | $k$ | | | | | | | | | |
|---|---|---|---|---|---|---|---|---|---|---|
| | 1 | 2 | 3 | 4 | 5 | 6 | 7 | 8 | 9 | 10 |
| 20 | 0.9756 | 0.8976 | 0.7636 | 0.5932 | 0.4168 | 0.2958 | 0.2718 | 0.2084 | 0.1915 | 0.1554 |
| 50 | 0.9956 | 0.9797 | 0.9473 | 0.8963 | 0.8268 | 0.7414 | 0.6447 | 0.5426 | 0.4411 | 0.3459 |
| 100 | 0.9989 | 0.9945 | 0.9852 | 0.9695 | 0.9465 | 0.9156 | 0.8766 | 0.8299 | 0.7762 | 0.7167 |
| 500 | 1.0000 | 0.9998 | 0.9993 | 0.9986 | 0.9975 | 0.9958 | 0.9937 | 0.9908 | 0.9873 | 0.9830 |

4.3. *Estimation.* In Section 3 the set of constraints (36) on the marginal likelihoods was derived for the case where all components are equivalent. These constraints can be used to perform a check of internal consistency of Markov chain Monte Carlo estimates of the marginal likelihoods $f(x|k)$, or of the marginal likelihoods implied by MCMC estimates of the posterior of $k$. The easiest way to check whether the constraints (36) are satisfied is to compute the $f_k^\dagger$ in (35) and see whether they are positive. As an example, Richardson and Green (1997) estimate a Bayesian mixture of univariate normals for the galaxy data set. They assume that all mixture components are equivalent, the prior on $k$ is $\pi(k) = k_{\max}^{-1}, k = 1, \ldots, k_{\max} = 30$, and the Dirichlet distributions on weights have hyperparameters $\alpha_{jk} = 1$. They report the reversible jump MCMC estimate of $\pi(k|x)$ contained in Table 5. Since the prior distribution of $k$ is uniform, the marginal likelihoods are proportional to the posterior of $k$. Substituting the above estimates of $\pi(k|x)$ for the $f_t$'s in (35), after disregarding the estimate for $k \geq 16$, produces, up to a proportionality constant, the $f_k^\dagger$'s implicit in Richardson and Green's estimate. These quantities are reported in Table 6. Three values of $\widehat{f_k^\dagger}$ are negative, for $k = 12, 13$ and 15. However, these violations are rather slight, almost within rounding error and occur for values of $k$ that account for little posterior probability and are, therefore, more difficult to estimate. Thus, if anything, the check gives support to Richardson and Green's estimate.

Table 4
*Bounds on $\pi(k|x)$ for several sample sizes $n$, $\pi(k) = k_{\max}^{-1}, k = 1, \ldots, k_{\max} = 50$, $\alpha = 0.5$*

| $n$ | $k$ | | | | | | | | | |
|---|---|---|---|---|---|---|---|---|---|---|
| | 1 | 2 | 3 | 4 | 5 | 6 | 7 | 8 | 9 | 10 |
| 20 | 0.7342 | 0.4684 | 0.2734 | 0.2575 | 0.1863 | 0.1783 | 0.1449 | 0.1343 | 0.1202 | 0.1030 |
| 50 | 0.8354 | 0.6477 | 0.4709 | 0.3229 | 0.2983 | 0.2618 | 0.2096 | 0.2047 | 0.1782 | 0.1664 |
| 100 | 0.8847 | 0.7456 | 0.6032 | 0.4703 | 0.3546 | 0.3166 | 0.2972 | 0.2610 | 0.2236 | 0.2189 |
| 500 | 0.9491 | 0.8833 | 0.8090 | 0.7306 | 0.6515 | 0.5742 | 0.5006 | 0.4320 | 0.3691 | 0.3392 |

TABLE 5
*Reversible jump MCMC estimate of $\pi(k|x)$ for the galaxy data set reported by Richardson and Green* (1997)

| $k$ | 1 | 2 | 3 | 4 | 5 | 6 | 7 | 8 |
|---|---|---|---|---|---|---|---|---|
| $\widehat{\pi(k|x)}$ | 0.000 | 0.000 | 0.061 | 0.128 | 0.182 | 0.199 | 0.160 | 0.109 |

| $k$ | 9 | 10 | 11 | 12 | 13 | 14 | 15 | $\geq 16$ |
|---|---|---|---|---|---|---|---|---|
| $\widehat{\pi(k|x)}$ | 0.071 | 0.040 | 0.023 | 0.013 | 0.006 | 0.003 | 0.002 | 0.003 |

Checking whether MCMC estimates of $f(x|k)$ or $\pi(k|x)$ satisfy the constraints only makes marginal use of the information supplied by them. This information can be more fully exploited by incorporating it in the estimation procedure. For instance, one could estimate the $f_k^{\dagger}$'s by MCMC methods and then use (26) to transform those estimates into estimates of the marginal likelihoods $f_k$'s. I will return to this point at the end of Section 4.4. Here I only sketch some approaches to transform estimates of the $f_k$'s into estimates that satisfy the inequalities (36).

TABLE 6
*Estimates, up to a proportionality constant, of $f_k^{\dagger}$ implicit in Richardson and Green* (1997) *MCMC estimate of $\pi(k|x)$, galaxy data set*

| $k$ | 1 | 2 | 3 | 4 | 5 | 6 | 7 | 8 |
|---|---|---|---|---|---|---|---|---|
| $\widehat{f_k^{\dagger}}$ | 0.0000 | 0.0000 | 0.0610 | 0.1194 | 0.1532 | 0.1413 | 0.0792 | 0.0352 |

| $k$ | 9 | 10 | 11 | 12 | 13 | 14 | 15 |
|---|---|---|---|---|---|---|---|
| $\widehat{f_k^{\dagger}}$ | 0.0167 | 0.0015 | 0.0035 | $-0.0005$ | $-0.0008$ | 0.0013 | $-0.0006$ |

TABLE 7
*Mode of* (43), *galaxy data,* $\widehat{\mathbf{f}}$ *is the Richardson and Green* (1997) *estimate given in Table* 5

| $k$ | 1 | 2 | 3 | 4 | 5 | 6 | 7 | 8 |
|---|---|---|---|---|---|---|---|---|
| $\mathbf{f}_k$ | 0.000 | 0.000 | 0.061 | 0.128 | 0.181 | 0.198 | 0.160 | 0.109 |

| $k$ | 9 | 10 | 11 | 12 | 13 | 14 | 15 |
|---|---|---|---|---|---|---|---|
| $\mathbf{f}_k$ | 0.071 | 0.041 | 0.023 | 0.013 | 0.007 | 0.003 | 0.002 |

Let $\mathbf{f} = (f(x|2), \ldots, f(x|k_{\max}))^{\top}$ be the vector of marginal likelihoods of the models with $k$ components, $k = 2, \ldots, k_{\max}$. Also, let $\widehat{\mathbf{f}}$ be the corresponding vector of MCMC estimates. When the mixture components parameters have conjugate prior distributions, $f_1 = f(x|1)$ can be computed exactly; if this is not the case, the vectors $\mathbf{f}$ and $\widehat{\mathbf{f}}$ also include $f(x|1)$ and its estimate. The estimates $\widehat{\mathbf{f}}$ might be directly available, as in the approaches of Nobile ([1994](#)), Carlin and Chib ([1995](#)), Raftery ([1996](#)) and Roeder and Wasserman ([1997](#)). Alternatively, they may only be computed up to a proportionality constant, from the prior on the number of components and an estimate $\widehat{\pi(k|x)}$ of its posterior, as in the approaches of Phillips and Smith ([1996](#)), Richardson and Green ([1997](#)) and Stephens ([2000](#)). In this latter case, the constraint proceeding from $\sum_{k=1}^{k_{\max}} \widehat{\pi(k|x)} = 1$ is disregarded. Estimates of the variability of $\widehat{\mathbf{f}}$ can be computed, either by replicating the MCMC runs or by using single run methods, such as batching and time series methods [see, e.g., Chapter 6 of Ripley ([1987](#)) or Geyer ([1992](#))]. It is assumed that as the MCMC sample size increases, the distribution of $\widehat{\mathbf{f}}$ approaches a multivariate normal

$$\widehat{\Sigma}^{-1/2}(\widehat{\mathbf{f}} - \mathbf{f}) \xrightarrow{\mathcal{D}} N(\mathbf{0}, I), \tag{42}$$

where $\widehat{\Sigma}$ is a consistent estimate of the variance–covariance matrix of $\widehat{\mathbf{f}}$. Let $R$ be the region where the constraints ([36](#)) are satisfied. If $\widehat{\mathbf{f}} \notin R$, an estimate of $\mathbf{f}$ which satisfies the constraints is the maximizer over $R$ of the likelihood $L(\mathbf{f})$ associated with ([42](#)). From a Bayesian viewpoint, this is equivalent to using $\widehat{\Sigma}$ as a plug-in estimate of $\Sigma$, employing $I_R(\mathbf{f})$ as the prior distribution of $\mathbf{f}$ and estimating $\mathbf{f}$ by the mode of its posterior distribution, which is proportional to

$$\exp\{-\tfrac{1}{2}(\mathbf{f} - \widehat{\mathbf{f}})^{\top}\widehat{\Sigma}^{-1}(\mathbf{f} - \widehat{\mathbf{f}})\}I_R(\mathbf{f}). \tag{43}$$

The posterior mode is the point in $R$ which is closest to $\widehat{\mathbf{f}}$ with respect to the metric induced by $\widehat{\Sigma}$. Hence, unless $\widehat{\mathbf{f}} \in R$, the mode will occur on the boundary of $R$, where the multivariate normal contours are tangent to $R$. The maximization of ([43](#)) is equivalent to the minimization of $(1/2)\mathbf{f}^{\top}\widehat{\Sigma}^{-1}\mathbf{f} - \widehat{\mathbf{f}}\widehat{\Sigma}^{-1}\mathbf{f}$ subject to $[\mathbf{b_2} \vdots \mathbf{b_3} \vdots \cdots \vdots \mathbf{b}_{k_{\max}}]\mathbf{f} \geq -\mathbf{b_1}f_1$, where the vector $\mathbf{b}_k$ has generic entry $b_{kt} = (-1)^{k+t}\binom{k}{t}a_{kt}I(k \geq t)$, $t = 2, \ldots, k_{\max}$. This is a simple problem in quadratic programming, for which software is publicly available; for instance, Goodall ([1995](#)) provides a basic S-PLUS implementation. Table [7](#) contains the $\mathbf{f}$ which maximizes ([43](#)) with $\widehat{\mathbf{f}}$ equal to the estimates of Richardson and Green ([1997](#)) given in Table [5](#).

Another estimate of $\mathbf{f}$, which satisfies the constraints ([36](#)) and does not lie on the boundary of $R$, is the mean of the distribution ([43](#)), which can be estimated by averaging independent draws from the posterior (43). However, drawing from the $N(\widehat{\mathbf{f}}, \widehat{\Sigma})$ distribution and using a rejection technique can be very inefficient, if $R$ is in the tail of the distribution. When this occurs, Gibbs sampling provides a more efficient alternative; working in terms of the distribution of the $f^\dagger$'s, a multivariate normal restricted to the positive orthant, leads to full conditional distributions that are univariate normals restricted to the positive reals. Table 8 contains an estimate of the posterior mean computed from 20,000 draws from (43), obtained using rejection, with $\widehat{\mathbf{f}}$ being Richardson and Green's (1997) estimate for the galaxy data. On the whole, the mean of (43) agrees with the estimate of Richardson and Green (1997), although it tends to give some more weight to models with a larger number of components. Table 9 displays the $f_k^\dagger$'s corresponding to the estimate of the mean of (43) given in Table 8. These estimates of the $f_k^\dagger$'s agree with those reported in Table 6 for values of $k$ up to 9, then they drop off much more regularly while remaining positive.

TABLE 8
*Estimate of the mean of* (43), *galaxy data,* $\widehat{\mathbf{f}}$ *is the Richardson and Green* (1997) *estimate given in Table* 5. *The estimate has been rescaled in order that* $\sum_k \mathbf{f}_k = 1$

| $k$ | 1 | 2 | 3 | 4 | 5 | 6 | 7 | 8 |
|---|---|---|---|---|---|---|---|---|
| $\mathbf{f}_k$ | 0.000 | 0.000 | 0.061 | 0.126 | 0.182 | 0.197 | 0.156 | 0.109 |

| $k$ | 9 | 10 | 11 | 12 | 13 | 14 | 15 | $\geq 16$ |
|---|---|---|---|---|---|---|---|---|
| $\mathbf{f}_k$ | 0.069 | 0.040 | 0.023 | 0.013 | 0.008 | 0.005 | 0.003 | 0.008 |

4.4. *The number of nonempty components.* Bayesian and classical analyses of the same data may lead to widely contrasting conclusions about the number of mixture components. A stylized account of a typical situation is as follows: a classical analysis identifies $\tilde{k}$ components as sufficient to provide a good fit to the data. On the other hand, the posterior of the number of components assigns considerable probability to values of $k > \tilde{k}$. Moreover, the posterior predictive distribution, conditional on $k$, of the next observation remains essentially the same for all $k \geq \tilde{k}$. Much of this divergence of conclusions derives from the use of the same term, in the two approaches, to denote different entities. In the Bayesian approach the parameter $k$ denotes the number of components in the mixture model, not the number of components from which data are actually observed. It is instead this second meaning that is attached to "number of components" in the classical approach; accordingly, determining the number of components amounts to finding $k$ such that $k$ mixture components afford a good fit of the data. The difference between the two approaches can be highlighted by positing a very small sample size, say $n = 3$; the classical approach will point at just one component, while the posterior of $k$ will be much the same as the prior. In the Bayesian approach it is quite possible for the posterior of $k$ to assign much probability to values larger than the number of components from which the data have originated. In fact, in Section 4.2 it was shown that, for a certain prior distribution, when $n = 82$ the posterior probability of three components is no larger than 0.8623, whatever the data are. This occurs because the posterior probabilities of four and more components cannot be too small, since they also account for allocation vectors with only three nonempty mixture components. As noted in Section 4.2, the strength of this link depends on the prior distribution of the mixture weights and it tends to abate as the sample size increases. However, the usefulness of the posterior of $k$, as a tool for selecting or estimating the number of components in a mixture, tends to be put in question by the fact that it may, to a very large extent, reflect probability mass associated with membership vectors that allocate observations to fewer than $k$ components.

TABLE 9
*Estimates of* $f_k^\dagger$ *corresponding to the mean of* (43) *given in Table* 8

| $k$ | 1 | 2 | 3 | 4 | 5 | 6 | 7 | 8 |
|---|---|---|---|---|---|---|---|---|
| $\widehat{f_k^\dagger}$ | 0.0000 | 0.0000 | 0.0612 | 0.1180 | 0.1536 | 0.1395 | 0.0766 | 0.0370 |

| $k$ | 9 | 10 | 11 | 12 | 13 | 14 | 15 |
|---|---|---|---|---|---|---|---|
| $\widehat{f_k^\dagger}$ | 0.0146 | 0.0033 | 0.0019 | 0.0007 | 0.0003 | 0.0002 | 0.0002 |

In summary, while the classical approach addresses the question:

Q1. *How many components are needed to fit the data well?*

The posterior of $k$ is suited to answer:

Q2. *How many components are likely to be in the model that generated the data?*

While Q2 is concerned with the number of components in the mixture, Q1 deals with the number of *nonempty* components. Since the Dirichlet prior on the mixture weights determines how likely empty components are to arise, it appears that the answer to Q2 depends on the prior specification more than the answer to Q1. This section seeks to pursue in a Bayesian way the objective of the classical approach, by deriving an expression for the posterior distribution of the number of nonempty components.

Let $h$ denote the number of nonempty components in the mixture. The joint prior distribution of the number of components $k$ and the membership vectors $g$ induces a prior on $h$. Since $h \leq k$, one has

$$f(h) = \sum_{k=h}^{\infty} \pi(k) f(h|k), \qquad h = 1, \ldots, n.$$

Let $\widetilde{\mathcal{G}}_h^k$ be the set of all membership vectors in $\mathcal{G}_k$ which assign observations to exactly $h$ components,

$$\widetilde{\mathcal{G}}_h^k = \bigcup_{t=h}^{k} \mathcal{G}_h^t. \tag{44}$$

Then the conditional distribution of $h$ given $k$ can be computed by summing $f(g|k, \alpha_k)$ over $\widetilde{\mathcal{G}}_h^k$,

$$f(h|k) = \sum_{g \in \widetilde{\mathcal{G}}_h^k} f(g|k, \alpha_k), \qquad h = 1, \ldots, k \wedge n. \tag{45}$$

The following proposition provides a representation of $f(h|k)$ which makes its computation feasible for sample sizes up to about 100; for larger samples sizes an estimate can be obtained by stochastic simulation.

PROPOSITION 4.2. *Consider the model of Section 2 under Assumptions* A.1 *and* A.2 *and suppose that all mixture components are equivalent. Let $d = d(n_1, \ldots, n_h)$ be the number of distinct entries in the vector $(n_1, \ldots, n_h)^\top$; also let $m_1, \ldots, m_d$ be the frequencies of the distinct $n_j$'s in $(n_1, \ldots, n_h)^\top$. Then*

$$
\begin{aligned}
f(h|k) = {} & \frac{\Gamma(k\alpha)}{\Gamma(k\alpha + n)} \binom{k}{h} \\
& \times \sum_{\substack{0 < n_1 \leq \cdots \leq n_h \\ n_1 + \cdots + n_h = n}} \binom{n}{n_1, \ldots, n_h} \binom{h}{m_1, \ldots, m_d} \\
& \times \prod_{j=1}^{h} \frac{\Gamma(\alpha + n_j)}{\Gamma(\alpha)}, \qquad h = 1, \ldots, k \wedge n.
\end{aligned}
\tag{46}
$$

Note that the sum in (46) does not involve $k$; this allows one to easily obtain $f(h|k)$ with $k > h$ from $f(h|k)$ with $k = h$. Therefore, one only needs to compute the sum in (46) at most $n$ times. The total number of terms in these $n$ sums is the number $p(n)$ of partitions of $n$ into integer summands without regard to order; tabulated values of $p(n)$ are in Table 24.5 of Abramowitz

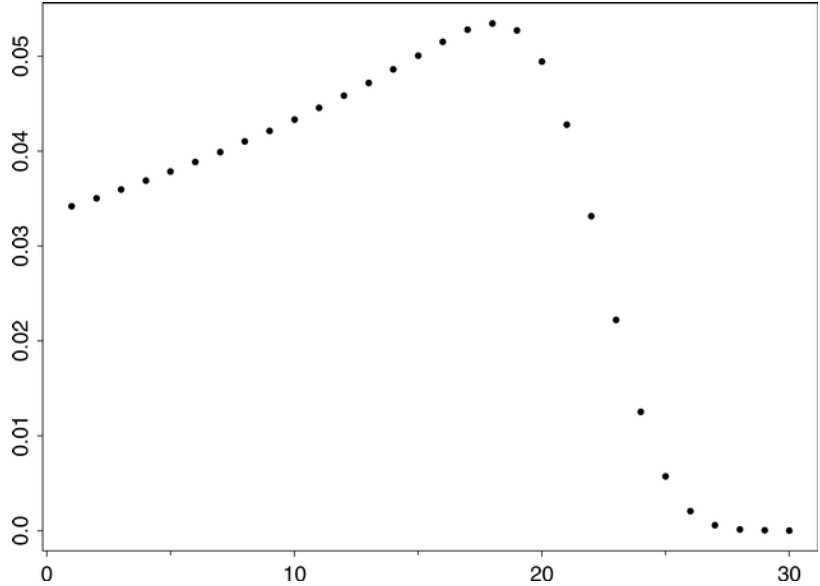

Fig. 1. *Prior distribution of the number $h$ of nonempty components when $n = 82$,*
$\pi(k) = 1/k_{\max}^{-1}, k = 1, \ldots, k_{\max} = 30$ *and* $\alpha = 1$.

and Stegun ([1964](#)). Figure [1](#) contains a plot of the prior distribution of $h$
corresponding to the prior used by Richardson and Green ([1997](#)) for the
galaxy data. The computation was done in Fortran and took six minutes on
a PC with a 1.1 GHz processor.

The posterior distribution of the number of nonempty components can
be written as

$$(47) \qquad f(h|x) = \sum_{k=h}^{\infty} \pi(k|x) f(h|k, x), \qquad h = 1, \ldots, n.$$

The following result provides a representation of the posterior of $h$ in terms
of the $f_h^\dagger$'s, the portions of the marginal likelihoods corresponding to no
empty components.

PROPOSITION 4.3.  *Consider the model of Section [2](#) under Assumptions
A.1 and A.2 and suppose that all mixture components are equivalent. Then*

$$(48) \quad f(h|x) = \frac{f_h^\dagger}{f(x)} \sum_{k=h}^{\infty} \pi(k) \binom{k}{h} \frac{\Gamma(k\alpha)}{\Gamma(k\alpha + n)} \frac{\Gamma(h\alpha + n)}{\Gamma(h\alpha)}, \qquad h = 1, \ldots, n.$$

Since the prior distribution of $h$ is only specified indirectly, through the
priors on $k$ and the mixture weights, one may prefer to consider, rather than the posterior of $h$, the marginal likelihood $f(x|h)$ for $h$ nonempty components. This quantity is readily derived from (48):

$$f(x|h) = \frac{f_h^\dagger}{f(h)} \sum_{k=h}^{\infty} \pi(k) \binom{k}{h} \frac{\Gamma(k\alpha)}{\Gamma(k\alpha + n)} \frac{\Gamma(h\alpha + n)}{\Gamma(h\alpha)}.$$

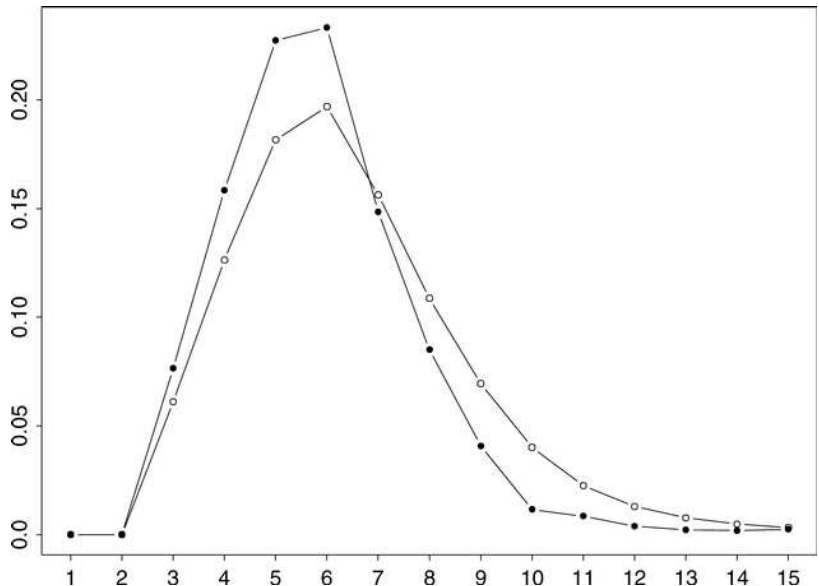

Fig. 2. *Estimates of $f(x|k)$ and $f(x|h)$ for the galaxy data, both normalized to sum to 1. Circles denote the estimate of $f(x|k)$ reported in Table 8; dots are the estimate of $f(x|h)$ obtained using the $f_h^\dagger$'s given in Table 9.*

Estimates of $f(x|h)$ are obtained by replacing the $f_h^\dagger$'s with the estimates produced in Section 4.3. Figure 2 displays estimates of $f(x|h)$, normalized to sum to 1, along with normalized estimates of the marginal likelihoods $f(x|k)$, for the galaxy data using the prior of Richardson and Green (1997). As one would expect, the marginal likelihoods of the number of nonempty components favor a smaller number of components than the posterior of $k$, effectively narrowing the plausible range of normal components in the observed data to between three and eight.

As a conclusion, note that the path here followed from estimates of the $f_k$'s to estimates of the $f_h^\dagger$'s to estimates of $f(x|h)$ can also be travelled in the opposite direction. For instance, it would be immediate to obtain estimates of $f(h|x)$ using Richardson and Green's (1997) reversible jump algorithm. These could then be turned, using (48), into estimates, up to a proportionality constant, of the $f_h^\dagger$'s and finally estimates of the marginal likelihoods $f_k$ automatically satisfying the constraints (36).

## APPENDIX: PROOFS

PROOF OF LEMMA 3.2.  The inverse image under $M_t$ of $\tilde{g} \in \mathcal{E}_h$ consists of all the $g \in \mathcal{G}_h^t$ which differ from $\tilde{g}$ only in that the nonempty components in each class can be any of the components in the class that are smaller than $t$, rather than being the first ones. If $m \neq i(t)$, there are $c(m,t)$ components in $\mathcal{C}(m)$ no larger than $t$, of which only $h_m$ are nonempty; there are $\binom{c(m,t)}{h_m}$ ways of selecting the nonempty components out of the $c(m,t)$ candidates. As $\mathcal{G}_h^t \subset \mathcal{G}_s^\star$, component $t$ is nonempty; this leaves $h_{i(t)} - 1$ nonempty components to be selected among $c(i(t),t) - 1$ candidates in $\mathcal{C}(i(t))$, yielding $\binom{c(i(t),t)-1}{h_{i(t)}-1}$ possible selections. Multiplying together the numbers of possible selections in the $N(t)$ classes yields (19).  □

PROOF OF LEMMA 3.3.  Use in (9) the partition of $\mathcal{G}_t^\star$ given in (17) to obtain $f_t^\star = \sum_{h \in \mathcal{H}_t} \sum_{g \in \mathcal{G}_h^t} f(x|t,g,\phi_t) f(g|t,\alpha_t)$. Replace the inner sum with the expression in (20): $f_t^\star = \sum_{h \in \mathcal{H}_t} \gamma_h^t \sum_{g \in \mathcal{E}_h} f(x|t,g,\phi_t) f(g|t,\alpha_t)$. Next recall that $\mathcal{E}_h \subset \mathcal{G}_s^\star$ and use C.1 and C.2: $f_t^\star = \sum_{h \in \mathcal{H}_t} a_{ts} \gamma_h^t \sum_{g \in \mathcal{E}_h} f(x|s,g,\phi_s) f(g|s,\alpha_s)$. Then use again (20) and then (21) to produce $f_t^\star = \sum_{h \in \mathcal{H}_t} a_{ts} (\gamma_h^t / \gamma_h^s) \sum_{g \in \mathcal{G}_h^s} f(x|s,g,\phi_s) f(g|s,\alpha_s) = \sum_{h \in \mathcal{H}_t} a_{ts} (\gamma_h^t / \gamma_h^s) f_h^\dagger$. From the partition of $\mathcal{H}_t$ in (15) it follows that $f_t^\star = \sum_{r=1}^t \sum_{h \in \mathcal{H}_r^t} a_{tr} (\gamma_h^t / \gamma_h^r) f_h^\dagger = \sum_{r=1}^t a_{tr} \sum_{h \in \mathcal{H}_r^r} I(h \in \mathcal{H}_r^t)(\gamma_h^t/\gamma_h^r) f_h^\dagger$, where the second equality uses the relationship in (16). Now (22) follows since, for all $h \in \mathcal{H}_r^r$, $I(h \in \mathcal{H}_r^t) = 0$ implies that $\gamma_h^t = 0$. To see this consider $h \in \mathcal{H}_r^r \backslash \mathcal{H}_r^t$. Since $s(h) = r$, $h \in \mathcal{H}_t$ would imply $h \in \mathcal{H}_r^t$ contrary to the hypothesis; hence $h \notin \mathcal{H}_t$ and from Definition 3.7 $\gamma_h^t = 0$.  □

PROOF OF THEOREM 3.2.  Substitute formula (22) in (11) to obtain $f_k = \sum_{t=1}^k a_{kt} \sum_{r=1}^t a_{tr} \sum_{h \in \mathcal{H}_r^r} (\gamma_h^t/\gamma_h^r) f_h^\dagger$. Next recall that $a_{kt} a_{tr} = a_{kr}$ and interchange the order of the two outer sums, $f_k = \sum_{r=1}^k a_{kr} \sum_{t=r}^k \sum_{h \in \mathcal{H}_r^r} (\gamma_h^t/\gamma_h^r) f_h^\dagger = \sum_{r=1}^k a_{kr} \sum_{h \in \mathcal{H}_r^r} f_h^\dagger (1/\gamma_h^r) \sum_{t=r}^k \gamma_h^t$. Finally, use (24) to produce (23). To prove (25) replace $f_k^\star$ in (12) with the expression provided by (22) with $t = k$.  □

PROOF OF COROLLARY 3.1.  Part (i) follows from Lemma 3.1 and Theorem 3.2. Equations (31) and (32) of part (ii) are obtained by replacing $a_{kh}$ in (26) and (27) with the expression given in (10) and using $\alpha_{0k} = k\alpha$. Part (iii) follows straightforwardly from part (ii) with $\alpha = 1$.  □

PROOF OF LEMMA 3.4.  Let $D = \{d_{kt}\}$ with $D = BC$. Then $D$ is lower triangular with generic element $d_{kt} = \sum_{r=1}^q b_{kr} c_{rt} = \sum_{r=1}^q (-1)^{r+t} b_{kr} b_{rt} = \sum_{r=t}^k (-1)^{r+t} b_{kr} b_{rt}$, with the last equality holding since $B$ is lower triangular. It then immediately follows that $D$ has unit diagonal elements, since so has $B$. Therefore, it only remains to show that $d_{kt} = 0$, $k > t$. Now use the definition of $b_{kt}$ and Condition C.2,

$$d_{kt} = \sum_{r=t}^{k} (-1)^{r+t} \binom{k}{r} a_{kr} \binom{r}{t} a_{rt}$$

$$= \sum_{r=t}^{k} (-1)^{r+t} \frac{k!}{r!(k-r)!} \frac{r!}{t!(r-t)!} a_{kt} = \binom{k}{t} a_{kt} \sum_{r=t}^{k} (-1)^{r+t} \binom{k-t}{k-r}.$$

Next, change the summation index to $j = r - t$ to obtain

$$d_{kt} = \binom{k}{t} a_{kt} \sum_{j=0}^{k-t} (-1)^{j+2t} \binom{k-t}{k-t-j} = \binom{k}{t} a_{kt} \sum_{j=0}^{k-t} (-1)^j \binom{k-t}{j} = 0,$$

as the sum $\sum_j$ is null because of a basic property of binomial coefficients [see, e.g., Abramowitz and Stegun (1964), page 10, Property 3.1.7].  □

PROOF OF PROPOSITION 3.4.   In the formula for $f_k$ given in (26) with $k > n$, replace $f_t^\dagger$ with the expression in (35) to produce

$$f_k = \sum_{t=1}^{n} \binom{k}{t} a_{kt} \sum_{r=1}^{t} (-1)^{t+r} \binom{t}{r} a_{tr} f_r$$

$$= \sum_{r=1}^{n} a_{kr} f_r \sum_{t=r}^{n} (-1)^{t+r} \binom{k}{t} \binom{t}{r}$$

(49)
$$= \sum_{r=1}^{n} a_{kr} f_r \binom{k}{r} \sum_{t=r}^{n} (-1)^{t+r} \binom{k-r}{k-t}.$$

Call $S$ the inner sum and rewrite it by changing the summation index to $j = t - r$ and making use of $\binom{k-r}{k-r-j} = \binom{k-r}{j}$,

(50)
$$S = \sum_{j=0}^{n-r} (-1)^j \binom{k-r}{j}.$$

Now, if $n-r$ is even, add $n-r-2j$ to the exponent of $(-1)$. This leaves $S$ unchanged, so that $S = \sum_{j=0}^{n-r} (-1)^{n-r-j} \binom{k-r}{j} = \binom{k-r-1}{n-r}$, where the last equality follows from a property of the binomial coefficients [see, e.g., Abramowitz and Stegun (1964), Section 24.1.1, Relations II.B]. If $n-r$ is odd, premultiply the sum $\sum_j$ in (50) by $-1$ and add $n-r-2j$ to the exponent of $(-1)$, yielding $S = -\binom{k-r-1}{n-r}$. Thus, in general, $S = (-1)^{n-r} \binom{k-r-1}{n-r}$. Finally, substituting the above expression of $S$ for the sum $\sum_t$ in the right-hand side of (49) and changing the index from $r$ to $t$ yields (37).  □

PROOF OF PROPOSITION 4.1. Rewrite (39) as follows:

$$\pi(k|x) = \sum_{h=1}^{k \wedge n} f_h^\dagger d_h \Big/ \sum_{h=1}^{n} f_h^\dagger b_h, \tag{51}$$

where $d_h = \pi(k)\binom{k}{h}a_{kh}$ and $b_h = \sum_{j=h}^{\infty}\pi(j)\binom{j}{h}a_{jh}$. It is immediate that a maximizer has $f_h^\dagger = 0$, $h > k \wedge n$, for otherwise $\pi(k|x)$ could be increased by simply setting these components to 0 and leaving the other ones unchanged. Suppose next that $\{f_h^\dagger\}_{h=1}^n$ has at least two nonzero components: there exist $t, r \in \{1, \ldots, k \wedge n\}$, $t \neq r$, such that $f_t^\dagger \neq 0$, $f_r^\dagger \neq 0$. Without loss of generality, assume that

$$\frac{b_r}{b_t} \geq \frac{d_r}{d_t}. \tag{52}$$

Define a new vector $\{\tilde{f}_h\}_{h=1}^n$ with $\tilde{f}_t = f_t^\dagger + (b_r/b_t)f_r^\dagger$, $\tilde{f}_r = 0$, $\tilde{f}_h = f_h^\dagger$, $h \neq t, r$. One can easily verify that replacing $f_h^\dagger$ with $\tilde{f}_h$ in the right-hand side of (51) leaves the denominator unchanged, while (52) ensures that the numerator does not decrease; $\sum_{h=1}^{k \wedge n} \tilde{f}_h d_h \geq \sum_{h=1}^{k \wedge n} f_h^\dagger d_h$. Therefore one can replace $f_h^\dagger$ with $\tilde{f}_h$ in (51), that is, select one of the nonzero components, set it to 0 and correspondingly adjust the other one, without decreasing $\pi(k|x)$. An appeal to induction completes the proof. □

PROOF OF PROPOSITION 4.2. Substituting in (45) $f(g|k, \alpha_k)$ from (3) and using the fact that all components are equivalent, one obtains

$$f(h|k) = \sum_{g \in \tilde{\mathcal{G}}_h^k} \frac{\Gamma(k\alpha)}{\Gamma(k\alpha + n)} \prod_{j=1}^{k} \frac{\Gamma(\alpha + n_j)}{\Gamma(\alpha)}.$$

The sum is over vectors $g$ with exactly $h$ nonempty components, so only $h$ terms in the products are not equal to 1. Since the terms in the sum do not depend on which components are nonempty, the sum is equal to $\binom{k}{h}$ times a sum over $\mathcal{G}_h^h$, the subset of $\mathcal{G}_h$ comprising vectors which allocate observations to all the $h$ mixture components. Therefore,

$$f(h|k) = \frac{\Gamma(k\alpha)}{\Gamma(k\alpha + n)} \binom{k}{h} \sum_{g \in \mathcal{G}_h^h} \prod_{j=1}^{h} \frac{\Gamma(\alpha + n_j)}{\Gamma(\alpha)}.$$

The terms in the above sum depend on $g$ only through $(n_1, \ldots, n_h)^\top$. Therefore one can replace the sum over $\mathcal{G}_h^h$ with a sum over all partitions of the $n$ observations in $h$ groups. Since to each partition $(n_1, \ldots, n_h)^\top$ there correspond $\binom{n}{n_1, \ldots, n_h}$ membership vectors in $\mathcal{G}_h^h$, one has

$$f(h|k) = \frac{\Gamma(k\alpha)}{\Gamma(k\alpha + n)} \binom{k}{h} \sum_{\substack{n_j > 0, \, j = 1, \ldots, h \\ n_1 + \cdots + n_h = n}} \binom{n}{n_1, \ldots, n_h} \prod_{j=1}^{h} \frac{\Gamma(\alpha + n_j)}{\Gamma(\alpha)}.$$

Finally, since the terms in the sum are invariant to a change in the order of the $n_j$'s, the sum above can be replaced by a sum over ordered $n_j$'s. As to each ordered vector $(n_1, \ldots, n_h)^\top$ there correspond $\binom{h}{m_1, \ldots, m_d}$ unordered ones, (46) follows.   $\square$

PROOF OF PROPOSITION 4.3.  The conditional distribution of $h$ given $k$ and $x$ in (47) can be obtained by summing the conditional distribution of $g$ given $k$ and $x$ over all membership vectors in $\mathcal{G}_k$ which allocate observations to exactly $h$ components; $f(h|k, x) = \sum_{g \in \widetilde{\mathcal{G}}_h^k} \{ f(x|k, g) f(g|k) \} / f(x|k)$. Substituting this expression in (47) produces

$$\begin{aligned}
f(h|x) &= \sum_{k=h}^{\infty} \frac{f(x|k)\pi(k)}{f(x)} \sum_{g \in \widetilde{\mathcal{G}}_h^k} \frac{f(x|k, g) f(g|k)}{f(x|k)} \\
&= \frac{1}{f(x)} \sum_{k=h}^{\infty} \pi(k) \sum_{g \in \widetilde{\mathcal{G}}_h^k} f(x|k, g) f(g|k).
\end{aligned}$$

(53)

Consider now the inner sum in (53):

$$\begin{aligned}
\sum_{g \in \widetilde{\mathcal{G}}_h^k} f(x|k, g, \phi_k) f(g|k, \alpha_k) &= \sum_{t=h}^{k} \sum_{g \in \mathcal{G}_h^t} f(x|k, g, \phi_k) f(g|k, \alpha_k) \\
&= \sum_{t=h}^{k} \sum_{g \in \mathcal{G}_h^t} f(x|t, g, \phi_t) f(g|t, \alpha_t) a_{kt} \\
&= \sum_{t=h}^{k} a_{kt} \gamma_h^t \sum_{g \in \mathcal{E}_h} f(x|t, g, \phi_t) f(g|t, \alpha_t),
\end{aligned}$$

where the first equality uses (44), the second one follows from Conditions C.1 and C.2 and the third uses (20). Now, when all components are equivalent $\mathcal{E}_h = \mathcal{G}_h^h$, so that using again Conditions C.1 and C.2 one obtains

$$\sum_{g \in \widetilde{\mathcal{G}}_h^k} f(x|k, g, \phi_k) f(g|k, \alpha_k) = \sum_{t=h}^{k} a_{kt} \gamma_h^t \sum_{g \in \mathcal{G}_h^h} f(x|h, g, \phi_h) f(g|h, \alpha_h) a_{th}$$

$$= a_{kh} f_h^\dagger \sum_{t=h}^{k} \gamma_h^t$$

with the second equality following from formula (28). Since $\gamma_h^t = \binom{t-1}{h-1}$ it follows that the inner sum in (53) equals $\binom{k}{h} a_{kh} f_h^\dagger$, so that

$$f(h|x) = \frac{f_h^\dagger}{f(x)} \sum_{k=h}^{\infty} \pi(k) \binom{k}{h} a_{kh}.$$

As an aside, note that the series in the right-hand side was already met in the denominator of (39). Substituting $a_{kh}$ with the expression in (10) and using $\alpha_{0k} = k\alpha$ yields (48). $\quad\square$

**Acknowledgments.** This paper develops an idea first presented in Section 4.4 of my Ph.D. dissertation. I thank Mark Schervish for being my thesis advisor and Larry Wasserman for repeated encouragement. Comments by Peter Green, Christian Robert, Mike Titterington and other participants at the Workshop in Statistical Mixtures and Latent-Structure Modelling held in Edinburgh, 28–30 March 2001, are also acknowledged. Helpful remarks by two anonymous referees and an Associate Editor have led to a sharper presentation.

## REFERENCES

ABRAMOWITZ, M. and STEGUN, I. A. (1964). *Handbook of Mathematical Functions.* National Bureau of Standards, Washington. MR167642

CARLIN, B. P. and CHIB, S. (1995). Bayesian model choice via Markov chain Monte Carlo methods. *J. Roy. Statist. Soc. Ser. B* **57** 473–484.

DIEBOLT, J. and ROBERT, C. P. (1994). Estimation of finite mixture distributions through Bayesian sampling. *J. Roy. Statist. Soc. Ser. B* **56** 363–375. MR1281940

GEYER, C. J. (1992). Practical Markov chain Monte Carlo (with discussion). *Statist. Sci.* **7** 473–503.

GOODALL, C. (1995). *S* functions for quadratic programming. Available from StatLib at http://lib. stat.cmu.edu/S/quadratic.

MENGERSEN, K. L. and ROBERT, C. P. (1996). Testing for mixtures: A Bayesian entropic approach. In *Bayesian Statistics 5* (J. M. Bernardo, J. O. Berger, A. P. Dawid and A. F. M. Smith, eds.) 255–276. Oxford Univ. Press. MR1425410

NOBILE, A. (1994). Bayesian analysis of finite mixture distributions. Ph.D. dissertation, Dept. Statistics, Carnegie Mellon Univ., Pittsburgh.

PHILLIPS, D. B. and SMITH, A. F. M. (1996). Bayesian model comparison via jump diffusions. In *Markov Chain Monte Carlo in Practice* (W. R. Gilks, S. Richardson and D. J. Spiegelhalter, eds.) 215–239. Chapman and Hall, London. MR1397970

RAFTERY, A. E. (1996). Hypothesis testing and model selection. In *Markov Chain Monte Carlo in Practice* (W. R. Gilks, S. Richardson and D. J. Spiegelhalter, eds.) 163–187. Chapman and Hall, London. MR1397966

RICHARDSON, S. and GREEN, P. J. (1997). On Bayesian analysis of mixtures with an unknown number of components (with discussion). *J. Roy. Statist. Soc. Ser. B* **59** 731–792. MR1483213

RIPLEY, B. D. (1987). *Stochastic Simulation.* Wiley, New York. MR875224

ROBERT, C. P. (1996). Mixtures of distributions: Inference and estimation. In *Markov Chain Monte Carlo in Practice* (W. R. Gilks, S. Richardson and D. J. Spiegelhalter, eds.) 441–464. Chapman and Hall, London. MR1397974

ROEDER, K. and WASSERMAN, L. (1997). Practical Bayesian density estimation using mixtures of normals. *J. Amer. Statist. Assoc.* **92** 894–902. MR1482121

STEPHENS, M. (2000). Bayesian analysis of mixture models with an unknown number of components: An alternative to reversible jump methods. *Ann. Statist.* **28** 40–74. MR1762903

TITTERINGTON, D. M., SMITH, A. F. M. and MAKOV, U. E. (1985). *Statistical Analysis of Finite Mixture Distributions.* Wiley, New York. MR838090

DEPARTMENT OF STATISTICS
UNIVERSITY OF GLASGOW
UNIVERSITY GARDENS
GLASGOW G12 8QW
UNITED KINGDOM
E-MAIL: agostino@stats.gla.ac.uk